\documentclass[12pt,reqno,a4paper,twoside]{article}

\usepackage{amsmath,amsthm,amstext,amscd,amssymb,euscript,mathrsfs}
\usepackage{calrsfs}
\usepackage{epsfig}

\newcommand{\Z}{\mathbb Z}
\newcommand{\R}{\mathbb R}
\newcommand{\N}{\mathbb N}

\newcommand{\E}{\mathbb E}

\renewcommand{\phi}{\varphi}

\newcommand{\tr}{\mbox{Tr}}

\newcommand{\pee}{\ensuremath{\mathbb{P}}}

\def\1{{\mathchoice {\rm 1\mskip-4mu l} {\rm 1\mskip-4mu l}
{\rm 1\mskip-4.5mu l} {\rm 1\mskip-5mu l}}}

\newtheorem{theorem}{{\small T}{\scriptsize HEOREM}}[section]
\newtheorem{corollary}{{\bf{\small C}{\scriptsize OROLLARY}}}[section]
\newtheorem{proposition}{{\bf{\small P}{\scriptsize ROPOSITION}}}[section]
\newtheorem{lemma}{{\bf{\small L}{\scriptsize EMMA}}}[section]
\newtheorem{remark}{{\bf{\small R}{\scriptsize EMARK}}}[section]
\newtheorem{definition}{{\bf{\small D}{\scriptsize EFINITION}}}[section]

\renewenvironment{proof}[1]
{\noindent{{\bf{\small{ P}{\scriptsize ROOF}}}.}\hspace{0.1cm} #1} {$\;\qed$\newline}

\newcommand{\beq}{\begin{eqnarray}}
\newcommand{\eeq}{\end{eqnarray}}

\newcommand{\ba}{\begin{align*}}
\newcommand{\ea}{\end{align*}}

\newcommand{\be}{\begin{equation}}
\newcommand{\ee}{\end{equation}}

\newcommand{\bl}{\begin{lemma}}
\newcommand{\el}{\end{lemma}}

\newcommand{\br}{\begin{remark}}
\newcommand{\er}{\end{remark}}

\newcommand{\bt}{\begin{theorem}}
\newcommand{\et}{\end{theorem}}

\newcommand{\bd}{\begin{definition}}
\newcommand{\ed}{\end{definition}}

\newcommand{\bp}{\begin{proposition}}
\newcommand{\ep}{\end{proposition}}

\newcommand{\bc}{\begin{corollary}}
\newcommand{\ec}{\end{corollary}}

\newcommand{\bpr}{\begin{proof}}
\newcommand{\epr}{\end{proof}}

\newcommand{\bi}{\begin{itemize}}
\newcommand{\ei}{\end{itemize}}

\newcommand{\ben}{\begin{enumerate}}
\newcommand{\een}{\end{enumerate}}

\newcommand{\caB}{{\mathcal B}}
\newcommand{\caC}{{\mathscr C}}

\newcommand{\caE}{{\mathrsfs E}}

\newcommand{\caF}{{\mathcal F}}

\newcommand{\caH}{{\mathcal H}}

\newcommand{\caL}{{\mathcal L}}
\newcommand{\caM}{{\mathcal M}}
\newcommand{\caN}{{\mathcal N}}

\newcommand{\caP}{{\mathcal P}}

\newcommand{\caR}{{\mathcal R}}
\newcommand{\caS}{{\mathcal S}}
\newcommand{\caT}{{\mathcal T}}

\begin{document}

\title{The abelian sandpile model on a random binary tree}

\author{
F.\ Redig $^{\textup{{\tiny(a)}}}$,
W. M. \ Ruszel $^{\textup{{\tiny(b)}}}$,
E.\ Saada $^{\textup{{\tiny(c)}}}$
\\
{\small $^{\textup{(a)}}$
Delft Institute of Applied Mathematics,}\\
{\small Technische Universiteit Delft}\\
{\small Mekelweg 4, 2628 CD Delft, Nederland}
\\
{\small $^{\textup{(b)}}$ Technical University Eindhoven}\\
{\small Department of Mathematics and Computer Science}\\
{\small P.O. Box 513,
5600 MB Eindhoven}\\
and \\
{\small Radboud University of Nijmegen}\\
{\small IMAPP}\\
{\small Heyendaalse weg 135,
6525 AJ Nijmegen }\\
{\small The Netherlands}\\
{\small $^{\textup{(c)}}$ CNRS, UMR 8145, Laboratoire MAP5,}\\
{\small
Universit\'{e} Paris Descartes, Sorbonne Paris Cit\'{e},} \\
{\small 45, Rue des Saints P\`{e}res, 75270 Paris Cedex 06, France}\\
}
\maketitle

\begin{abstract}
We study the abelian sandpile model on a random binary tree.
Using a transfer matrix approach introduced by  Dhar \& Majumdar,
we prove exponential decay of correlations,
and in a small supercritical region  (i.e., where the branching process survives with positive probability)  exponential decay of avalanche sizes.
This shows a phase transition phenomenon between exponential decay
and power law decay of avalanche sizes.
Our main technical tools are:
(1) A recursion
for the ratio between the numbers of weakly and strongly allowed
configurations which is proved to have a well-defined stochastic solution;
(2) quenched and annealed estimates of the
eigenvalues of a product of $n$ random transfer matrices.
\end{abstract}
\section{Introduction}
The abelian sandpile model (ASM) is a thoroughly studied model both in the physics
and in the mathematics literature see e.g. \cite{dharo,MRSindien,
jaraio,redigo} for recent review papers on the subject. In physics, it serves as a paradigmatic model
of self-organized criticality (SOC). SOC is usually referred to as the phenomenon
that the model exhibits power law decay of correlations or avalanche sizes, without
fine-tuning
any external parameters such as temperature or magnetic field.
In mathematics, the ASM is connected to several combinatorial objects such as spanning
trees, graph-orientations, dimers, and it has an interesting abelian group structure.

The ASM has been studied on the Bethe lattice (i.e., the rootless binary tree) in \cite{dhar}. Via a recursive
analysis, based on a transfer matrix method, the authors in \cite{dhar}
arrive at exact expressions of various quantities of interest, such as
the single height distribution, correlation functions of height variables,
and avalanche size distribution.

There are various motivations to consider the ASM on {\em random graphs}.
As an example, we mention integrate-and-fire models in neuroscience (see e.g.\ \cite{lev}
and related papers),
where the connections between neurons are updated after a neuron has fired.
The typical connection structure of a network of firing neurons is therefore
generically not translation invariant, and time dependent.
As a first approximation, one can quench the randomness of the connection
graph and study the firing of neurons on the derived random graph.
So far, the ASM model has been studied on small world graphs from the
physicist's perspective, using a renormalization group approach
\cite{smallworld}.
In the mathematics literature, there are recent studies on so-called
``cactus'' graphs \cite{cactus}.

In this paper, we start this study of the ASM on random graphs
with the ASM on a random tree, for the sake
of simplicity chosen to be a realization of a binary branching process with branching probability $p\in[0,1]$.
We use the transfer-matrix method of \cite{dhar} to express
relevant quantities such as the correlation of height variables and
the avalanche size distribution in terms of the eigenvalues of
 an ad hoc  product of \textit{random} matrices. This is the fundamental difference
between the Bethe lattice case and
the random tree, namely the fact that the transfer matrices depend randomly
on the vertices and instead of having to deal with the $n$-th power of a simple
two by two matrix, one has to control the product of $n$
random matrices.

The crucial quantity entering the transfer matrices is the so-called
characteristic ratio, which is the ratio between the numbers of weakly
and strongly allowed configurations. This ratio is equal to $1$ in the infinite Bethe
lattice  for every vertex and it is close to $1$ for vertices belonging
to a finite subset of the Bethe lattice which are far away from the ``boundary''
 (see later on for precise statements).
In our case, we
show that for an infinite random tree
the characteristic ratio is a well-defined random variable, uniquely
determined by a stochastic recursion. The transfer matrices will then
contain elements with that distribution.
We also consider
deterministic trees that are strict subsets of
the binary tree where the
characteristic ratio can be computed explicitly.
Next, we prove the exponential decay of correlation of height variables
(as in the Bethe lattice case), and show that for a branching probability $p$
sufficiently small, but still supercritical ($p>1/2$),  i.e., the branching process survives with positive probability,  avalanche
sizes decay exponentially. This shows a transition
between exponential decay of avalanche sizes, for $p$ small, and
power law decay for $p$ close to (possibly only equal to) one.

Our paper is organized as follows.
First we recall some basic material about the ASM on trees and the
recursive technique developed in \cite{dhar}.
Second we study the recursion for the characteristic ratio
and show it has a unique solution for the random binary tree.
Finally we give quenched and annealed estimates of the eigenvalues of the product of $n$ random
matrices, which we apply
in the study of correlation of height variables
and avalanche sizes.
\section{Abelian sandpile model on subtrees of the full binary tree}\label{sec:ASM-on-trees}
We summarize here the basic and standard objects of the abelian sandpile model
on a (general) tree.
More details can be found e.g.\ in \cite{dhar, MRS}.
\subsection{Rooted and unrooted random trees}\label{rootedandnot+random}
We denote by $\caB_n$ the rooted binary tree of $n$ generations,
and by $\caB_\infty$ the rooted infinite binary tree.
For a more general tree $\caT$ we write $\caT^i$ if we want to indicate that the tree has root $i$.
The rootless infinite binary tree or Bethe lattice is then obtained by joining two
infinite rooted binary trees by a single edge connecting their roots.

A random binary tree of $N$ generations with branching probability $p\in [0,1]$ is a random subset  $\caT_N$  of $\caB_N$ obtained as follows.
Starting from the root, we add two new vertices, each connected with a single edge to the root (resp.\ no vertices), with probability $p$  (resp.\ $(1-p)$),  and we iterate this
from every new vertex independently for $N$ generations. By letting $N\to\infty$ we obtain the full binary branching process.
Joining two independent infinite copies of this process by a single edge connecting their roots creates the
rootless random binary tree.
 This last procedure is of course identical to create the non-random rootless binary tree from non-random rooted binary trees.
\subsection{Height configurations and legal topplings}\label{subsec:height-legal}
For $\caT$ a finite subtree of the Bethe lattice, height configurations
on $\caT$ are elements $\eta\in\{1,2,\ldots\}^\caT:=\caH_\caT$. For $\eta\in\caH_\caT$ and
$u\in\caT$,
$\eta_u$ denotes the height at vertex $u$.
A height configuration $\eta\in\caH_\caT$
is \textit{stable} if $\eta_u\in \{1,2,3\}$ for all $u\in\caT$.
Stable configurations are collected in the set
$\Omega_\caT=\{1,2,3\}^{\caT}$.

For a configuration $\eta\in\caH_\caT$, we define the \textit{toppling
operator} $T_u$ via
\[
\left(T_u (\eta)\right)_v=\eta_v-\Delta_{uv}
\]
where $\Delta$ is the \textit{toppling matrix}, indexed by vertices $u,v\in\caT$ and
defined by
\be\label{top}
\Delta_{uu}=3, \Delta_{uv}=-1\ \mbox{if}\ u,v\
\mbox{are neighbors in}\ \caT
\ee
($u,v$ neighbors in $\caT$ means that an edge of $\caT$ connects $u$ to $v$).
In words, in a toppling at $u$, 3 grains are removed from $u$, and
every neighbor of $u$ receives one grain.

A toppling at $u\in\caT$ is called \textit{legal} if $\eta_u>3$.
A sequence of legal topplings is a composition $T_{u_n}\circ\ldots\circ T_{u_1} (\eta)$ such that
for all $k=1,\ldots,n$ the toppling at $u_k$ is legal
in $T_{u_{k-1}}\circ\ldots\circ T_{u_1} (\eta)$. The \textit{stabilization}
of a configuration $\eta\in\caH_\caT$ is defined
as the unique stable configuration $\caS(\eta)\in\Omega_\caT$
that arises from $\eta$ by a sequence of legal topplings.
\subsection{Addition operator and Markovian dynamics}
For $\caT$ a finite subtree of the Bethe lattice, and for $u\in\caT$,
the \textit{addition operator} is the map $a_u:\Omega_\caT\to\Omega_\caT$
defined via
\be\label{ax}
a_u \eta = \caS (\eta+\delta_u)
\ee
where $\delta_u\in\{0,1\}^\caT$ is such that  $\delta_u(u)=1$ and $\delta_u(z)=0$ for $z\in\caT,z\not=u$.
In other words, $a_u\eta$ is the effect of an addition of a single grain at $u$ in $\eta$,
followed by stabilization.

The addition operators commute, i.e., $a_ua_v=a_va_u$.
This is the well-known and crucial \textit{abelian property} of the sandpile model.

The dynamics of the sandpile model is then
the discrete-time Markov chain $\{\eta(n), n\in \N\}$
on $\Omega_\caT$ defined via
\be\label{markov}
\eta(n)= \prod_{i=1}^n a_{X_i} \eta(0)
\ee
where $X_i$ are i.i.d.\ uniformly chosen vertices of $\caT$.

Given a stable height configuration $\eta$
and $u\in\caT$, we define
the \textit{avalanche} $Av(u,\eta)$ induced by addition at $u$ in $\eta$ to be
the set of vertices in $\caT$ that have to be toppled in the course of the stabilization
of $\eta+\delta_u$.
\subsection{Recurrent configurations and stationary measure}\label{subsec:rec-stat}
The \textit{recurrent} configurations of the sandpile model
form a subset of the stable
configurations defined as follows. A configuration $\eta\in\Omega_\caT$
 contains a \textit{forbidden subconfiguration} (FSC)
if there exists a subset $S\subset\caT$ such that for all
$u\in S$, the height in $u$ is less than or equal to
the number of neighbors of $u$ in $S$. The restriction of $\eta$ to $S$ is then called a FSC. A configuration is \textit{allowed}
if and only if it does not contain a FSC. Recurrent configurations coincide
with allowed ones, and
are collected in the set $\caR_\caT$.

 We denote by $\caP(\Omega_\caT)$ the set of probability measures on $\Omega_\caT$.
The Markov chain \eqref{markov} has a unique stationary probability measure
$\mu_\caT\in\caP(\Omega_\caT)$  which is
the uniform measure on the set $\caR_\caT$
\[
\mu_\caT= \frac{1}{|\caR_\caT|}\sum_{\eta\in\caR_\caT} \delta_\eta
\]
where $\delta_\eta$ is the point mass concentrated on the configuration
$\eta$.
\subsection{Specific properties of the sandpile model on a tree}\label{treesection}
In this section, for the sake of self-containdness, we briefly summarize some basic facts from the paper
\cite{dhar} which we need later on.
In a subset of the Bethe lattice, the \textit{distance between two vertices} is defined as the length
(i.e., the number of edges) of the shortest path
joining them.
A vertex is a \textit{surface vertex} if it has a number of neighbors strictly less than
$3$.
\subsubsection{Weakly and strongly allowed subconfigurations}\label{weakstrong}
The class of allowed configurations can be divided into \textbf{weakly} and \textbf{strongly} allowed ones.
Let $T$ be a rooted finite tree with root $u$ and extend it with
one vertex $v$ and one edge $<u,v>$. Consider an allowed  configuration $\xi$ on $T$. We put
height 1 at vertex $v$  and investigate the derived configuration $\xi^{\prime}$ on $T \cup \lbrace v \rbrace$ such that $\xi$
is the restriction of $\xi'$ to $T$, and $\xi'_v=1$.
If  $\xi^{\prime}$ has no FSC, then we call $\xi$ \textit{strongly} allowed on $T$, otherwise \textit{weakly} allowed.
 We give in Figure 1 an example of weakly and strongly allowed configurations.
\begin{figure}[ht]
\begin{center}
\includegraphics[scale=1]{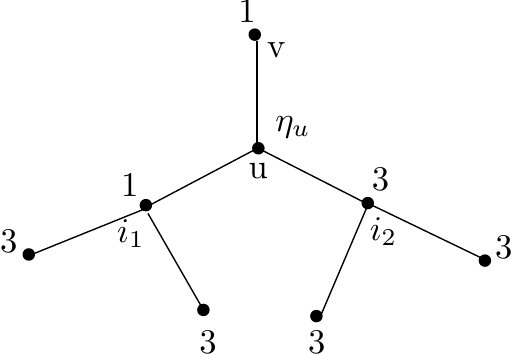}
\caption{Example of a tree with root $u$}
\end{center}
\end{figure}
On one hand, if $\eta_u=2$, then there exists a forbidden subconfiguration on $S=\lbrace v, u, i_1\rbrace$.
On the other hand, if $\eta_u=3$, there are no forbidden subconfigurations.
\subsubsection{Characteristic ratio and recursion}\label{cr+rec}
Let $\caT$ be a finite tree, rooted or not.
A key quantity in the analysis of \cite{dhar} is
\textit{the characteristic ratio} $x(\caT)\in [1/2, 1]$ of $\caT$
between weakly and strongly allowed configurations. For the empty tree $\caT=\varnothing$ we put $x(\varnothing)=0$.

For an infinite tree $\caT$ we say that the characteristic ratio is well-defined if the limit
$\lim_{\caT'\uparrow\caT} x(\caT')$ exists, where $\lim_{\caT'\uparrow\caT}$ is taken
along the net of finite subtrees. That is, $\lim_{\caT'\uparrow\caT} x(\caT')=a$ means that for every
$\varepsilon>0$, there exists $\caT'_0\subset\caT$ a finite subtree of $\caT$ such that, for all
finite subtrees $\caT'\supset\caT'_0$, we have $|x(\caT')-a|<\varepsilon$.

The characteristic ratio satisfies a recursion property both for rooted and unrooted trees:
\begin{equation}\label{dharrec}
 x(T) = \frac{(1 + x(T^1))(1 + x(T^2))}{2 + x(T^1) + x(T^2)}
\end{equation}
where $T^1,T^2$ are two non-intersecting subtrees of the tree $T$ defined as follows.
 If $T=T^{o}$ is a rooted tree with root $o$,
then $T^1,T^2$ are obtained by
deleting the root $o$ and splitting the tree into two subtrees whose roots are the two descendants of the root $o$,
see Figure 2.
\begin{figure}[htu]\label{WR}
\begin{center}
 \includegraphics[scale=0.55]{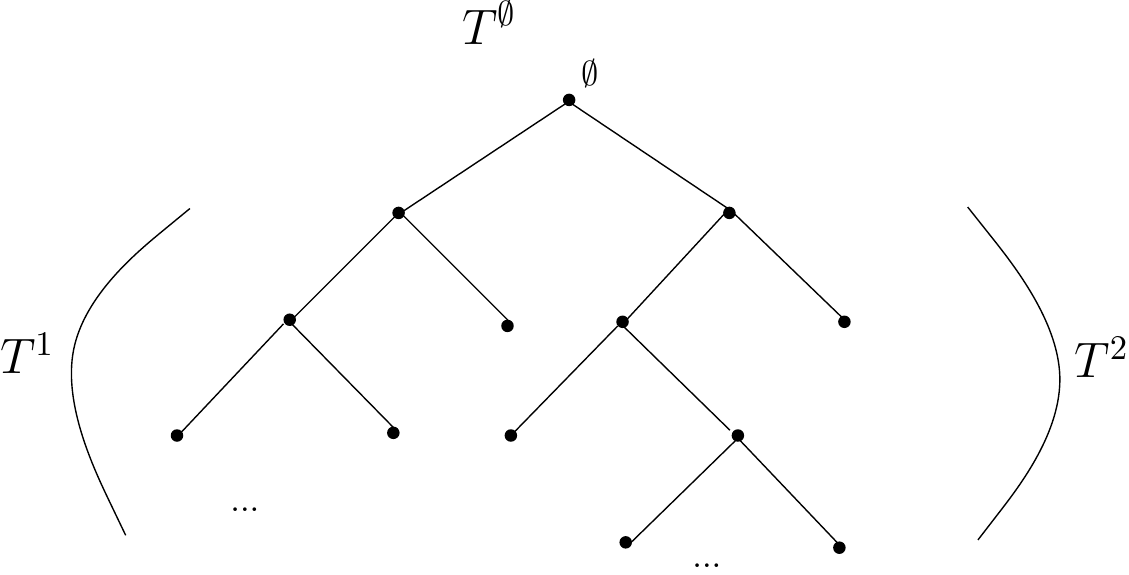}
\caption{Example of a tree with root $o$ which splits into $T^1$ and $T^2$}
\end{center}
\end{figure}
If $T$ is an unrooted tree, we pick one of its edges, $<u,v>$,
we delete this edge and split $T$ into two rooted subtrees $T^1$ and $T^2$ whose roots are $u$ and $v$, see Figure 3.
\begin{figure}[htu]\label{NR}
\begin{center}
 \includegraphics[scale=0.55]{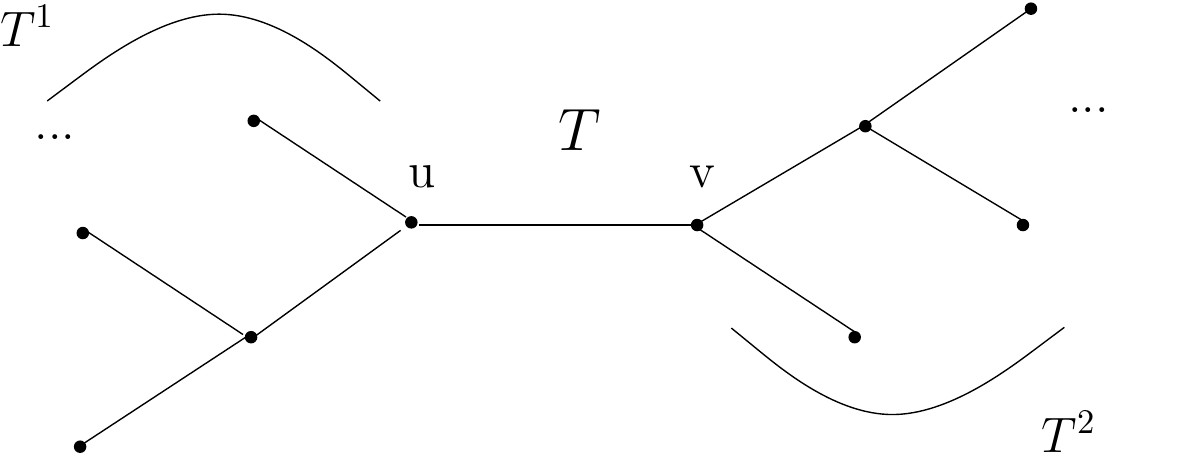}
\caption{Example of a rootless tree which splits into $T^1$ and $T^2$}
\end{center}
\end{figure}
This recursion \eqref{dharrec} holds for finite trees, and by passing to the limit for infinite trees, provided
the limits defining $x(T),x(T^1),x(T^2)$ exist.
 Dhar \& Majumdar obtain from %the recursion 
  \eqref{dharrec} that when the distance from the root of $T$ to the nearest surface vertex tends to infinity (what they call ``deep in the lattice''), the characteristic ratio $x(T)$ tends to 1.
The authors then derive several explicit quantities (such as height probabilities of a vertex
$u$
deep in the lattice) by replacing the characteristic ratio by $1$.
\subsubsection{Transfer matrix approach}\label{transmatapp}
 The transfer matrix approach allows to compute the two-point correlation functions and
avalanche size distribution. Let $u,v$ be two vertices
in the tree at mutual distance $n$. They determine $n+3$ subtrees $T_1,\ldots T_{n+3}$,
see Figure 4. The number of allowed configurations when fixing the heights at
$u$ and $v$ can then be obtained via a product $\caM$ of $n+3$
two by two matrices
with elements
determined by the characteristic ratios of the subtrees
$T_1,\ldots, T_{n+3}$ (see Section 4 below for the precise form of these matrices).
\begin{figure}[ht]\label{TM}
\begin{center}
\includegraphics[scale=0.8]{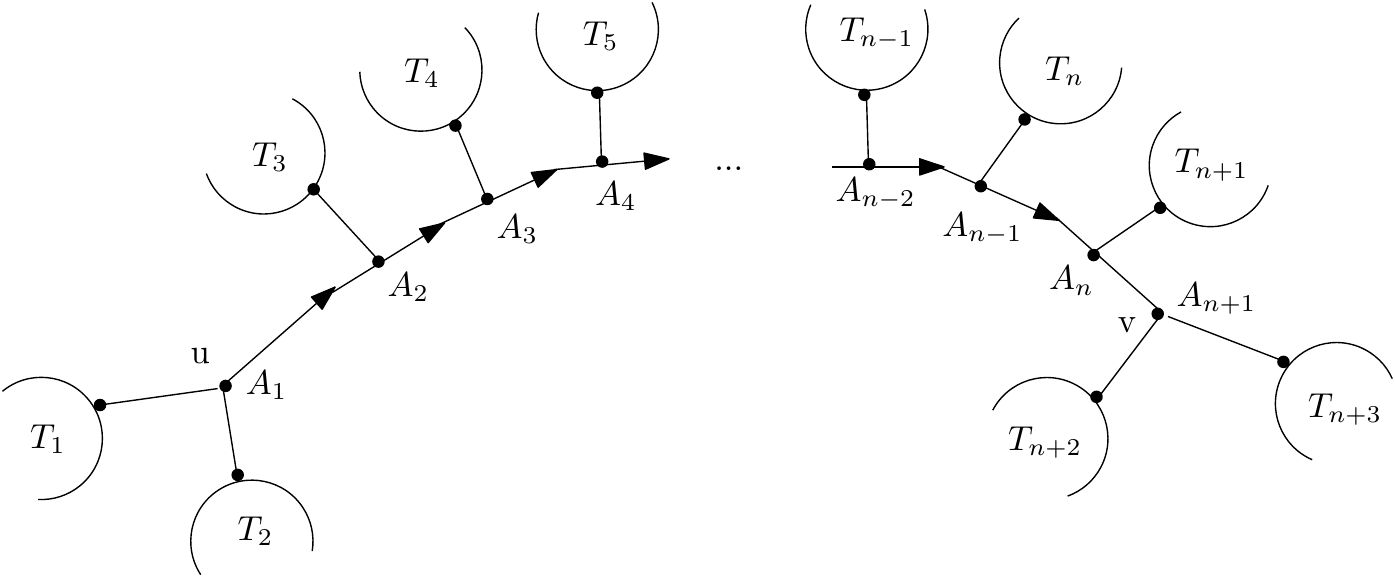}
\caption{Example of a tree with $n+3$ subtrees}
\end{center}
\end{figure}
In particular for the Bethe lattice, all the matrices involved in this product are equal
to
$
\left(\begin{array}{cc}
2&2\\
1& 3
\end{array}\right)
$
because the characteristic ratios of all the involved subtrees are equal to one.
Let $\caT$ be a finite or infinite subtree of the full binary tree. The two-point correlation function, i.e., the probability that two vertices $u,v$ at mutual distance $n$ have height $i$ resp. $j$, for $i,j\in\{1,2,3\}$, is equal to (see \cite[Section 5, eq. (5.11)]{dhar})
\be
\mu_{\caT}(\eta_u=i,\eta_v=j) = \mu_{\caT}(\eta_u=i)\mu_{\caT}(\eta_v=j) + a_{i,j}\frac{\lambda_{-}(\caM)}{\lambda_{+}(\caM)}\label{covbinary}
\ee
where $\lambda_{-}(\caM)/\lambda_{+}(\caM)$ is the ratio of the smallest and largest eigenvalues of the matrix $\caM$ and $a_{i,j}$ are some numerical constants depending on $i,j$. For the binary tree,  $\lambda_{-}(\caM)=1$ and $\lambda_{+}(\caM)=4^{n}$.

The avalanche size distribution is determined by the inverse of $\lambda_{+}(\caM)$. For the binary tree,
upon addition of a grain at a vertex $u$, the probability that the avalanche
$Av(u,\eta)$ is a given connected subset $\caC$ of cardinality $n$ containing $u$
is equal to
\be
\mu_{\caT}(Av(u,\eta)=\caC) = \frac{C}{4^{n}},\label{avdis}
\ee
for some constant $C$, independent of the shape of the subset $\caC$.
Since there are $4^n n^{-3/2} (1+ o(1))$  connected subsets of the Bethe lattice of cardinality $n$ containing $u$, one concludes
that for large $n$ (see  \cite[Section 6, eq. (6.13), (6.14)]{dhar}),
\be
\mu_{\caT}(|Av(u,\eta)|=n) \approx n^{-3/2}
\ee
i.e., the tail of the avalanche size distribution decays like $n^{-3/2}$.
 \section{Some characteristic ratios}\label{charsec}

On the Bethe lattice $x(\caT^i)=1$ for the (infinite) rooted subtrees ($i=1,2,3$) attached
to every vertex since every $\caT^i$ is an infinite
rooted binary tree (see Subsection \ref{cr+rec}). This property considerably simplifies the analysis
of \cite{dhar} and is no longer valid in the inhomogeneous or random 
cases studied here.
\subsection{The characteristic ratio of the random binary tree}\label{crt}
Let $\caT_n$ be the random  binary tree of $n$ generations starting
from a single individual (the ``zero-th'' generation) at time $n=0$
%, where every individual has two children with probability $p$ and
%no children with probability $(1-p)$, $p\in [0,1]$
 (cf. Section \ref{rootedandnot+random}).
Then for $n>0$, $x(\caT_n)$ satisfies the recursive identity
\cite[Section 3, eq. (3.12)]{dhar}
\be\label{recur}
x(\caT_n)=
f(x(\caT_{n-1}^1),x(\caT_{n-1}^2))
\ee
where
 \be\label{functionre}
f(u,v)= \frac{(1+u)(1+v)}{2+u+v}= \left(\frac{1}{1+u}+\frac{1}{1+v}\right)^{-1}
\ee
and
$\caT^i_{n-1}$ are the (possibly empty) subtrees
emerging from the (possibly absent) individuals of the first generation.
For a tree $\caT_0$ consisting of a single point we have
\be\label{cr1}x(\caT_0)= \frac{1}{2}\ee
 because  heights $2,3$ are strongly allowed
and  height $1$ is weakly allowed (cf. Section \ref{weakstrong}).
This value $1/2$ can also be obtained from the recursion \eqref{recur}
by viewing a single point
as connected to two empty subtrees (for which $x(o)=0$).
Notice that if $u,v\in [0,1]$, then $f(u,v)\leq 1$. Therefore, we view $f$ as
a function from $[0,1]^2$ onto $[1/2,1]$.

\bl\label{lemma f}
For every finite subtree $\caT$ of the Bethe lattice,
the ratio $x(\caT)\in [1/2,1]$.
Moreover, on $[0,1]^2$ the function $f$ defined by \eqref{functionre}
is symmetric, i.e., $f(u,v)=f(v,u)$, and increasing in $u$ and $v$, i.e., for $u_1 \leq u_2$,
\begin{equation}
f(u_1,v) \leq f(u_2,v), \text{ for all } v
\end{equation}
and analogously in the other argument.
\el
\bpr
The proof is straightforward and left to the reader.
\epr
\bp\label{contractionbinary}
There exists a random variable $X_\infty$ such that
$x(\caT_n)\to X_\infty$ in distribution as $n\to\infty$.
\ep
\bpr
Denote by $\mu_n\in\caP(\Omega_{\caT_n})$ the distribution of $x(\caT_n)$.
The recursion \eqref{recur}
induces a corresponding recursion on the distributions $\mu_{n+1}=\caF(\mu_n)$.
We show that $\caF$ is a contraction on the set $\caP([1/2,1])$ of probability measures on $[1/2,1]$
 endowed with the Wasserstein distance. This
implies that it has a unique fixed point $\mu^*$ and from every initial $\mu_0$,
$\mu_n\to\mu^*$ in Wasserstein distance and thus weakly.

Let $g$ be a Lipschitz function. Then we have
\begin{eqnarray} \label{rec exp}
%\mathbb{E}[g(x(\caT_n))]
\int g(x(\caT_n)) \mu_n(d(x(\caT_n)))& =& \sum_{T }g(x(T))%\mathbb{P}
\mu_n(x(\caT_n) = x(T))\nonumber\\
& =& \sum_{i\in \{0,1,2\}}g(x(T^i)) %\mathbb{P}
\mu_n(x(\caT_n) = x(T^i)) 
\nonumber \\
&+& \sum_{ T \neq T^0,T^1,T^2} g(x(T))%\mathbb{P}
\mu_n(x(\caT_n) = x(T))
\end{eqnarray}
where $T^0,T^1,T^2$ are the three trees which
cannot be split into two subtrees both non-empty, see Figure 5.
\begin{figure}[ht]\label{T0T1T2}
\begin{center}
 \includegraphics[scale=0.8]{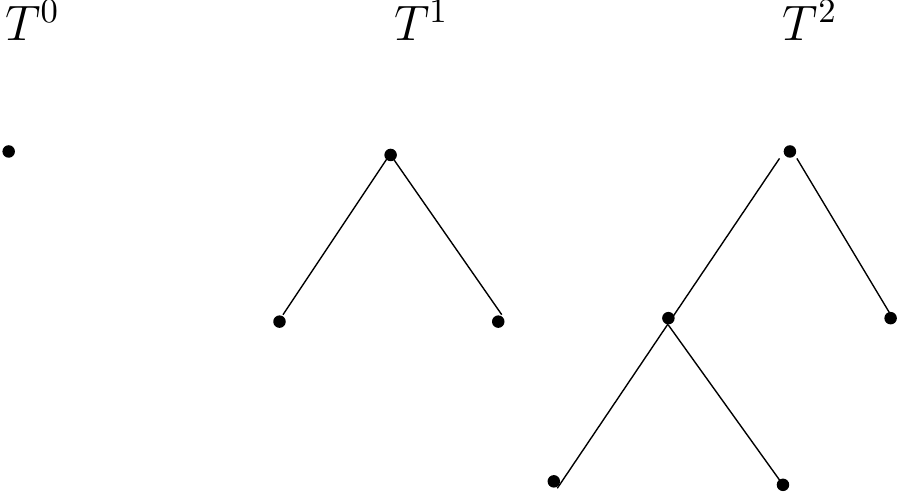}
 \caption{The trees $T^0,T^1,T^2$}
\end{center}
\end{figure}
Every other tree $T$ (appearing in the last sum of \eqref{rec exp}) can be split into two subtrees
 both not reduced to a single point.
Then using \eqref{recur}, the expression in \eqref{rec exp}  becomes
\begin{eqnarray}
%\mathbb{E}[g(x(\caT_n))]
&&\int g(x(\caT_n)) \mu_n(d(x(\caT_n)))\nonumber\\
& =& g \biggl(\frac{1}{2} \biggr) (1-p) + g \circ f \biggl(\frac{1}{2},\frac{1}{2}  \biggr) p(1-p)^2
\nonumber\\
 &+& 2p^2(1-p) %\mathbb{E}
\int \biggl(g \circ f \biggl(x(\caT_{n-1}),\frac{1}{2}\biggr)\biggr)\mu_{n-1}(d(x(\caT_{n-1})))
 \nonumber\\
 &+& p^3 %\mathbb{E}
 \int \biggl (g\circ f(x(\caT_{n-1}),x(\caT_{n-1})) \biggr ) \mu_{n-1}(d(x(\caT_{n-1}))).
\end{eqnarray}
\bl\label{lemme3.2}
The function $\caF$ on $\caP([1/2,1])$ defined by
\begin{eqnarray} \label{operator F}
\int g(x) \mathcal{F}(\mu)(dx)  &=&
g \biggl(\frac{1}{2} \biggr) (1-p) + g\biggl(\frac{3}{4} \biggr) p(1-p)^2
\nonumber\\
& +& 2p^2(1-p)
\int g\circ f\biggl(x,\frac{1}{2} \biggr)\mu(dx)\nonumber\\
& +& p^3 \int \int g \circ f(x,y)\mu(dx)\mu(dy)
\end{eqnarray}
is a contraction on $\caP([1/2,1])$
endowed with the Wasserstein distance. The contraction factor is bounded
from above by $8/9$.
\el
\bpr
Denote by $\caL$ the set of Lipschitz functions $g:[1/2,1]\to\R$ with Lipschitz constant less than or equal to one, i.e.,
such that $|g(x)-g(y)|\leq |x-y|$ for all $x,y$.
We use the following two formulas for the Wasserstein distance of two elements $\mu,\nu$ of $\caP([1/2,1])$ \cite{dudley}:
\be\label{was1}
d(\mu,\nu) =\sup \left\{\left|\int gd\mu - \int gd\nu\right|: g\in\caL\right\},
\ee
\be\label{was2}
d(\mu,\nu)= \inf\left\{\int |x-y| \pee(dx dy): \pee_1=\mu, \pee_2=\nu\right\}
\ee
where in the last right hand site the infimum is over all couplings $\pee$ with first marginal $\pee_1$ (resp.\ second marginal $\pee_2$)
equal to $\mu$ (resp.\ $\nu$).

To estimate $d(\caF (\mu),\caF(\nu))$ we start with the first formula \eqref{was1}
\begin{equation}\label{diff F}
 \begin{split}
& d(\mathcal{F}(\mu),\mathcal{F}(\nu))=
  \sup_{g\in\caL}\biggl | \int g(x)\mathcal{F}(\mu)(dx) - \int g(y)\mathcal{F}(\nu)(dy) \biggr |  \\
& \\
& \leq 2p^2(1-p)\sup_{g\in\caL}\biggl | \int g\circ f\biggl(x,\frac{1}{2}\biggr)\mu(dx) - \int g\circ f\biggl(y,\frac{1}{2} \biggr)\nu(dy) \biggr |\\
& \\
& + p^3 \sup_{g\in\caL}\biggl | \int \int g \circ f(x,x^{\prime}) \mu(dx)\mu(dx^{\prime}) - \int \int g \circ f(y,y^{\prime}) \nu(dy)\nu(dy^{\prime}) \biggr |\\
& \\
& =: 2p^2(1-p)\sup_{g\in\caL} A(\mu,\nu,g) + p^3\sup_{g\in\caL} B(\mu,\nu,g).
 \end{split}
\end{equation}
Now use the definition of $f$ and the fact that $x,y\in [1/2,1]$ to estimate
\[
\left|f\left(x,\frac12\right)-f\left(y,\frac12\right)\right|\leq \frac14 |x-y|.
\]
This gives, using the Lipschitz property of $g$ and  a coupling $\pee$ of $\mu$ and $\nu$: 
\begin{eqnarray}
A(\mu,\nu,g)
&=&
\biggl | \int g\circ f\biggl(x,\frac{1}{2}\biggr)\mu(dx) - \int g\circ f\biggl(y,\frac{1}{2} \biggr)\nu(dy) \biggr |
\nonumber\\
&\leq&
 \int \biggl| g\circ f\biggl(x,\frac{1}{2}\biggr) -  g\circ f\biggl(y,\frac{1}{2} \biggr)\biggr| \pee(dx, dy)
\nonumber\\
&\leq &
 \int \biggl| f\biggl(x,\frac{1}{2}\biggr) -   f\biggl(y,\frac{1}{2} \biggr)\biggr| \pee(dx,dx^{\prime},dy,dy^{\prime})
\nonumber\\
&\leq &
\frac14\int |x-y| \pee(dx, dy).
\label{Amunug}
\end{eqnarray}
Taking now the infimum over all couplings $\pee$, using \eqref{was2} to bound \eqref{Amunug},
and taking the supremum over $g\in \caL$ yields
\be\label{brol1}
\sup_{g\in \caL} A(\mu,\nu,g)\leq \frac14 d(\mu,\nu).
\ee
To estimate the term $B(\mu,\nu,g)$, use the elementary bound (since $x,x',y,y'\in [1/2,1]$):
\[
|f(x,x')-f(y,y')|\leq \frac19(4|x-y| + 4 |x'-y'|).
\]
We then have, using the Lipschitz property of $g$, and a coupling $\pee'= \pee'(dxdx';dydy') $ of $\mu\otimes\mu$ and $\nu\otimes \nu$ 
\begin{eqnarray*}
&&\biggl | \int \int g \circ f(x,x^{\prime}) \mu(dx)\mu(dx^{\prime}) - \int \int g \circ f(y,y^{\prime}) \nu(dy)\nu(dy^{\prime}) \biggr |
\\
&=&
\int\int |f(x,x')-f(y,y')| \pee'(dxdx';dydy')
\\
&\leq & \frac49\int|x-y| \pee'(dxdx';dydy') + \frac49\int |x'-y'| \pee'(dxdx';dydy').
\end{eqnarray*}
Taking the supremum over $g\in\caL$, and infimum over the couplings $\pee'$, we find
\be\label{brol2}
\sup_{g\in\caL} B(\mu,\nu,g) \leq \frac89 d(\mu,\nu).
\ee
Combining \eqref{brol1}, \eqref{brol2} with \eqref{diff F} we arrive at
\[%be\label{cont11}
d(\caF (\mu),\caF(\nu)) \leq \left(\frac12p^2 (1-p) + \frac89 p^3\right)d(\mu,\nu)\leq \frac{8}{9} d(\mu,\nu)
\]%ee
where in the final inequality we used the elementary bound
\[
\frac12p^2 (1-p) + \frac89 p^3= p^2 \left(\frac12 + \frac{7p}{18}\right)\leq \frac{8}{9}.
\]
\epr
 This proof of Lemma \ref{lemme3.2} completes the proof of Proposition \ref{contractionbinary}.
\epr
\br\label{rk:intro-binomial}
The random {\rm binomial} tree is such that every vertex has two children with probability $p^2$, 1 child with probability $2p(1-p)$ and 0 children with probability $(1-p)^2$, $p \in [0,1]$. For this case the same idea as in Lemma \ref{lemme3.2} gives a recursion leading to a contraction in the Wasserstein distance and hence $x(\caT_n)$ has also a unique limiting distribution.
\er
\subsection{The characteristic ratio of some deterministic trees}\label{crdettrees}
The recursion \eqref{recur} allows also to compute the characteristic
ratio for certain (deterministic) infinite subsets of the full binary
tree in terms
of iterations of sections of $f$.
\subsubsection{Infinite branch}\label{infbranch}
First, consider a single branch of length $n$, i.e.,
the tree $\caT^{min}_n$ consisting of a root and $n\geq 1$ generations of two individuals each
 (see Figure 6).
\begin{figure}[ht]
\begin{center}
 \includegraphics[scale=0.8]{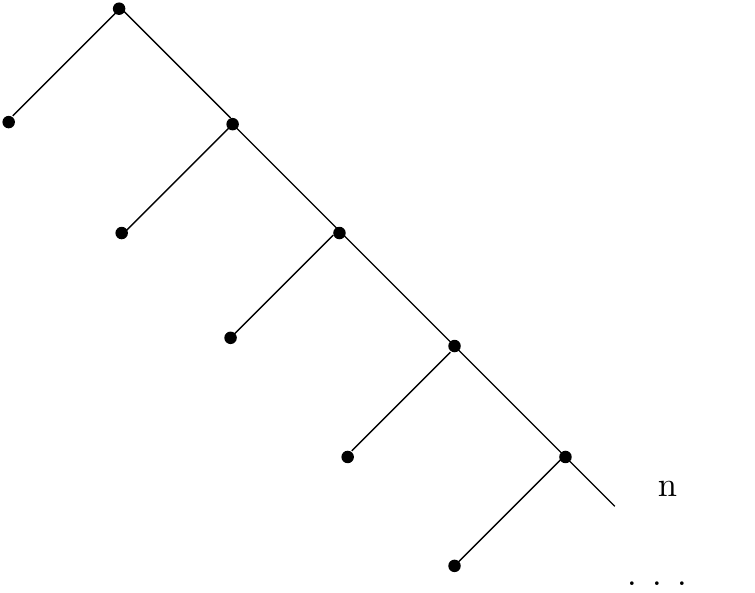}
\caption{Example of a tree with a single branch of $n$ generations}
\end{center}
\end{figure}
Using \eqref{cr1} we obtain the recursion
\be\label{recursbac}
x\left(\caT^{min}_{n+1}\right)= f\left(\frac12, \caT^{min}_n\right)
\ee
which by Lemma \ref{lemma f} gives that $x\left(\caT^{min}_{n}\right)$
is monotonically increasing in $n$ with limit
\be\label{xinfty}
x_{\infty}=x\left(\caT^{min}_{\infty}\right)= \frac{-1+\sqrt{7}}{2}
\ee
which is the unique positive solution of
\be\label{xfixed}
x= f\left(\frac12, x\right).
\ee
\br
We can generalize this ``backbone'' tree to a general ``backbone-like''
tree where at each point
the same finite tree $T$ is attached (it is a singleton in the backbone case).
The characteristic ratio equals the positive solution of the fixed-point equation
\be
x^*= f(x(T), x^*)
\ee
 where $x(T)$ denotes the characteristic ratio of the tree $T$, i.e.,
\be
 x^* = \frac{1}{2} (-1 + \sqrt{5 + 4 x(T)}).
\ee
\er
\subsubsection{Finite perturbations of a single branch}
Attaching a finite subtree $T$ of the full binary tree $\caB_\infty$
at level $n$ in the infinite single branch $\caT^{min}_\infty$ leads
to a tree $\caT^{per}_n$ with characteristic ratio
\be\label{perrec}
x(\caT^{per}_n)=\phi\circ\phi\ldots\circ\phi (f(x(T),x(\caT^{min}_\infty)))
\ee
where $\phi(x)= f(1/2, x)$ (cf. \eqref{xfixed})  is applied
$n$ times.
This shows that inserting a finite tree $T$ at level $n$ has an effect
on the characteristic ratio that
vanishes in the limit $n\to\infty$, exponentially fast in $n$.

Moreover, since $x(T)\geq 1/2$ and $x\mapsto f(x,y)$ is monotone
for all $y$ (by Lemma \ref{lemma f}), we
have from \eqref{perrec}
\be\label{mono}
x(\caT^{per}_n)\geq \phi\circ\phi\ldots\circ\phi \left(f\left(\frac12,x(\caT^{min}_\infty)\right)\right)=x(\caT^{min}_\infty).
\ee
{}From \eqref{mono}  and \eqref{xinfty} we conclude that for every {\em infinite} subtree
$\caT_\infty
\subset\caB_\infty$
for which $x(\caT_\infty)$ exists,
\[
\frac{-1+\sqrt{7}}{2}=x(\caT^{min}_\infty)\leq x(\caT_\infty)\leq x(\caB_\infty)=1.
\]
\section{Transfer matrix and eigenvalues: uniform estimates}\label{tramateig}
In the analysis of the two point correlation function and of
the avalanche size distribution, one is confronted with
the problem of estimating the minimal and maximal eigenvalues (denoted by $\lambda_{-}(\caM)$ and \ $\lambda_+(\caM)$) of a
product of matrices of the form
\be\label{matprod}
\caM(x_1,\ldots,x_n)=\prod_{i=1}^n M(x_i)
\ee
with
\be
M(x_i)=\left(\begin{array}{cc}
1+x_i&1+x_i\\
1&2+x_i\\
\end{array}\right)
\ee
where the $x_i$'s are the characteristic ratios of some recursively defined subtrees.

When $x_i=1$ for all $i$, this exactly corresponds to the analysis
in \cite{dhar} (see the above Subsection \ref{transmatapp}). More precisely, in that case, one needs the minimal and maximal
eigenvalues (denoted by $\lambda_{-}$ and \ $\lambda_+$) of
\[
\caM(1,1\ldots, 1) =\left(\begin{array}{cc}
2&2\\
1&3\\
\end{array}\right)^n
\]
which are $\lambda_{-}=1$ and $\lambda_+=4^n$. This
leads to a decay of the covariance proportional to
$1/4^n= \lambda_{-}/\lambda_{+}$ and decay of avalanche size
asymptotically proportional to
\[
\left(\frac{1}{4^n}\right)\frac{4^n}{n^{3/2}}= \frac{1}{\lambda_{+}} A_n
\]
where
\be\label{animal}
A_n=|\{ \caC\subset\caB_\infty: \caC\ \mbox{connected},\ |\caC|=n, \ \caC\ni o \}|
\ee
denotes the cardinality of the set of connected clusters of size $n$ containing 
the origin $o$. In the general case the decay of covariance can be estimated
by the ratio $\lambda_{-}(\caM)/\lambda_+(\caM)$, and for the decay of the avalanche
size distribution, one needs to estimate $\lambda_+(\caM)$ as
well as the analogue of $A_n$.
In the case
of a branching process, the $x_i$'s appearing in the matrix
$\caM(x_1,\ldots,x_n)$ are independent random variables
with distribution $\mu^*$ defined in the proof of Proposition \ref{contractionbinary}.

In this section
we therefore concentrate on the estimation of the eigenvalues
of a matrix of the form $\caM(x_1,\ldots,x_n)$ for general
$x_i$'s.
\bl\label{eigenlem}
\begin{enumerate}
\item For all $n$ and all $x_i\in [1/2, 1]$, ($1\le i\le n$), the eigenvalues of
$\caM(x_1,\ldots,x_n)$ are non-negative.
\item We have the inequality
\be\label{2}
\frac{\lambda_-(\caM)}{\lambda_+(\caM)}\leq C \left(\frac{4}{9}\right)^n.
\ee
\item We have
\be\label{3}
\frac{\lambda_-(\caM)}{\lambda_+(\caM)}=\frac{ \det (\caM(x_1,\ldots,x_n))}{\tr(\caM(x_1,\ldots,x_n))^2} (1+o(1))
\ee
\end{enumerate}
where $o(1)$ tends to zero as $n\to\infty$, uniformly
in the choice of the $x_i$'s.
\el
\bpr
The eigenvalues are given by
\be\label{eigen}
\lambda_{\pm}(\caM)=\frac{1}{2}\left( a\pm\sqrt{a^2-4b}\right)
\ee
with $a= \tr(\caM)$ and $b=\det (\caM)$
(we abbreviate $\caM$ for $\caM(x_1,\ldots,x_n)$
in these expressions).
For $n=1$, $a^2-4b= 5+4x\geq 0$. Hence
\be\label{singleeigen}
\lambda_{\pm}= \frac{3+2x\pm\sqrt{5+4x}}{2}.
\ee
For $n\geq 2$, we estimate,
as in \cite{MRS}
\[
a=\tr(\caM)\geq \prod_{i=1}^n(2+x_i)
\]
and using
\be\label{detM}
b=\lambda_{+}(\caM)\lambda_{-}(\caM)=\prod_{i=1}^n (1+x_i)^2
\ee
we have
\be\label{1}
\frac{b}{a^2}\leq \prod_{i=1}^n \left(1+\frac{1}{1+x_i}\right)^{-2}\leq \left(\frac{4}{9}\right)^n.
\ee
In particular, $4b\leq a^2$ which implies that the eigenvalues are real
and non-negative ($a\geq 0$). Inequality \eqref{2}
then follows from \eqref{1}, \eqref{eigen}.

Given that (by \eqref{1}) $h=b/a^2$ tends to zero as $n\to\infty$ at a speed
at least $C (4/9)^n$, we have
\[
\frac{\lambda_-(\caM)}{\lambda_+(\caM)}=\frac{1-\sqrt{1-4h}}{1+\sqrt{1+4h}}= h(1+O(h))
\]
which proves the third statement.
\epr

Lemma \ref{eigenlem}
shows that the ratio $\lambda_-(\caM)/\lambda_+(\caM)$ behaves
in leading order as $\det(\caM)/(\tr(\caM)^2)$.
To estimate this ratio, we need to estimate $\tr(\caM)$ from below.
We start with a useful representation of $\caM (x_1,\ldots,x_n)$.
\bl
Define\label{formtrace}
\be\label{e1}
E_1=
\left(\begin{array}{cc}
1&1\\
0&1
\end{array}
\right), \quad
E_2=
\left(\begin{array}{cc}
0&0\\
1&1
\end{array}
\right)
\ee
Then we have
\be\label{form}
\caM(x_1,\ldots,x_n)= \sum_{\alpha=(\alpha_1,\ldots,\alpha_n)\in \{0,1\}^n} \prod_{i=1}^n y_i^{\alpha_i}\caE(\alpha)
\ee
where $y_i= (1+x_i)$ and
where
\be\label{eea}
\caE(\alpha)=\prod_{i=1}^n\left( E_1^{\alpha_i}E_2^{1-\alpha_i}\right).
\ee
\el
\bpr
We have
\[
\caM (x_1,\ldots,x_n)= \prod_{i=1}^n (y_i E_1+ E_2).
\]
The result then follows from expansion of this product.
\epr
\bl\label{lowerlem}
\ben
\item For all $n\geq 1$,
\be\label{prodma}
E_1^n=
\left(\begin{array}{cc}
1&n\\
0&1
\end{array}
\right),\quad
E_2^n =E_2, \quad
E_1^n E_2=
\left(\begin{array}{cc}
n&n\\
1&1
\end{array}
\right).
\ee
\item For all $r\geq 1$, $k_1,\ldots,k_r,k_{r+1}\geq 0$,
\be\label{case1}
\tr\left(\prod_{i=1}^r (E_1^{k_i}E_2)\right)= \prod_{i=1}^r (1+k_i)
\ee
and
\be\label{case2}
\tr\left(\left(\prod_{i=1}^r (E_1^{k_i}E_2)\right) E_1^{k_{r+1}}\right)=\left(\prod_{i=2}^r (1+k_i)\right) (1+k_1+k_{r+1}).
\ee
\een
\el
\bpr
Identity \eqref{case2} follows
from \eqref{case1} and invariance of the trace
under cyclic permutations.
To prove \eqref{case1}, use the expression in \eqref{prodma} for $E_1^nE_2$, and
estimate the diagonal elements of the product
\[
\left(\begin{array}{cc}
k_1&k_1\\
1&1
\end{array}
\right)
\left(\begin{array}{cc}
k_2&k_2\\
1&1
\end{array}
\right)
\ldots
\left(\begin{array}{cc}
k_r&k_r\\
1&1
\end{array}
\right)=
\left(\prod_{i=2}^r (1+k_i)\right)
\left(\begin{array}{cc}
k_1&k_1\\
1&1 \end{array}
\right)
\]
which implies the result.
\epr
\bp\label{propo}
For all $x_1,\ldots,x_n\in [1/2,1]$
we have 
\be\label{lowertrace}
\tr (\caM)\geq \left(\prod_{i=1}^n (1+2y_i)\right) 2^{-n(16/25)}.
\ee
As a consequence we have the following uniform
upper bound
\be\label{upperratio}
\frac{\lambda_-(\caM)}{\lambda_+(\caM)}\leq C\gamma^n
\ee
where
\[
\gamma=\left(\frac{4}{25} 2^{32/25}\right)\approx 0.38854
\]
\ep
\bpr
Remember definition \eqref{eea} of $\caE(\alpha)$ for $\alpha=(\alpha_1,\ldots,\alpha_n)\in \{0,1\}^n$. Using
Lemma \ref{lowerlem} we see that
\[
\tr(\caE(\alpha))\geq 2^{\caN(\alpha)}
\]
where
$\caN(\alpha)=\sum_{i=1}^{n-1} \alpha_i(1-\alpha_{i+1})$, - with the
convention $\alpha_{n+1}=1$ - is the
number of intervals of successive $1$'s in the configuration $\alpha$.
Hence by \eqref{form} we obtain the lower bound 
\beq\nonumber
\tr(\caM(x_1,\ldots,x_n))&=& \sum_{\alpha=(\alpha_1,\ldots,\alpha_n)\in \{0,1\}^n} \left(\prod_{i=1}^n y_i^{\alpha_i}\right)\tr(\caE(\alpha))\\
&\geq& \sum_{\alpha=(\alpha_1,\ldots,\alpha_n)\in \{0,1\}^n} \left(\prod_{i=1}^n y_i^{\alpha_i}\right)2^{\sum_{i=1}^n \alpha_i(1-\alpha_{i+1})}.\label{tussen1}
\eeq
Next since  $\sum_{\alpha=(\alpha_1,\ldots,\alpha_n)\in \{0,1\}^n}\left(\prod_{i=1}^n(2y_i)^{\alpha_i}\right)=\prod_{i=1}^n(1+2y_i)=:Z$, we rewrite
\be\label{Zbound}
\sum_{\alpha=(\alpha_1,\ldots,\alpha_n)\in \{0,1\}^n} \left(\prod_{i=1}^n y_i^{\alpha_i}\right)2^{\sum_{i=1}^n \alpha_i(1-\alpha_{i+1})}
=Z \E_a \left(2^{-\sum_{i=1}^n\alpha_i\alpha_{i+1}}\right)
\ee
 where
\[
\E_a (\psi(\alpha))=\frac{1}{Z} \sum_{\alpha=(\alpha_1,\ldots,\alpha_n)\in \{0,1\}^n}\left(\prod_{i=1}^n(2y_i)^{\alpha_i} \right)\psi(\alpha)
\]
defines a $(y_i)_i$ dependent probability measure on the $\alpha$'s.
Now apply Jensen's inequality
in \eqref{Zbound} to obtain
\begin{eqnarray*}
&&\frac1{Z}\sum_{\alpha=(\alpha_1,\ldots,\alpha_n)\in \{0,1\}^n} \left(\prod_{i=1}^n y_i^{\alpha_i}\right)2^{\sum_{i=1}^n \alpha_i(1-\alpha_{i+1})}\\
&\geq & 2^{-\E_a(\sum_{i=1}^n\alpha_i\alpha_{i+1})}= 2^{-\sum_{i=1}^n\frac{4y_iy_{i+1}}{(1+2y_i)(1+2y_{i+1})}}.
\end{eqnarray*}
 Finally, use $y_i=1+x_i\in [3/2, 2]$ to estimate
\[
\sum_{i=1}^n\frac{4y_iy_{i+1}}{(1+2y_i)(1+2y_{i+1})}\leq \frac{16n}{25}
\]
which gives the following uniform lower bound for $\tr(\caM)$:
\[
\frac{1}{Z}\tr(\caM(x_1,\ldots,x_n))\geq 2^{-n(16/25)}.
\]
Combining this with $\det(\caM)=\prod_{i=1}^n a^2_i$ we
obtain
\[
\frac{\lambda_-(\caM)}{\lambda_+(\caM)}\leq (1+o(1))2^{n(32/25)}\prod_{i=1}^n\left(\frac{y_i}{1+2y_i}\right)^2.
\]
Finally we use $y_i=(1+x_i)$ with $x_i\in [1/2,1]$ and the fact that
$x\mapsto (1+x)(3+2x)^{-1}$ is
increasing to estimate
\[
\left(\frac{y_i}{1+2y_i}\right)^2\leq \frac4{25}
\]
which implies \eqref{upperratio}.
\epr

\section{Transfer matrix: annealed estimates}
In this section we look at the eigenvalues of $\caM(x_1,\ldots,x_n)$
where now the $x_i$'s are i.i.d. with a law
$\mu$ on $[1/2,1]$. We denote by $\pee$ the joint law of the
$x_i$'s and by $\E$ the corresponding expectation.

We start with the following lemma.
\bl
For all $\gamma\geq 0$ the eigenvalues
of
\[
A(\gamma)=
\left(
\begin{array}{cc}
\gamma &\gamma\\
1 & 1+\gamma
\end{array}
\right)
\]
given by
\be\label{consteig}
\Lambda_{\pm}(\gamma)= \frac{2\gamma+1\pm \sqrt{4\gamma+1}}{2}
\ee
are non-negative.
\el
\bpr
Elementary computation.
\epr
\bt
Let $\lambda_{\pm}(\caM)$ denote the largest, resp.\ smallest,
eigenvalues of $\caM (x_1,\ldots,x_n)$ where the $x_i$'s are 
i.i.d.\ with a law supported on $[1/2,1]$. Denote
\be
Y_n= \frac1{n} \log\left(\frac{\lambda_+(\caM)}{\lambda_-(\caM)}\right).
\ee
Then we have:
\ben
\item Concentration property: there exists $C>0$ such that for all $\varepsilon >0$
\be\label{conc}
\pee\left(\left(Y_n-\E(Y_n)\right)>\varepsilon\right)\leq e^{-C\varepsilon^2 n}.
\ee
\item The limits
\be\label{osel}
L_{\pm}= \lim_{n\to\infty}\frac1n\log\lambda_{\pm}(\caM)
\ee
 exist and satisfy
$L_{\pm}=\E(L_{\pm})$ almost surely; moreover
\be\label{prodeig}
L_++L_-= 2\E\left(\log (1+x_1)\right).
\ee
\item Upper bound
\be\label{annealup}
\lim_{n\to\infty} Y_n=L_+-L_- =\lim_{n\to\infty} \E(Y_n)\leq \log\frac{\Lambda_+ (\gamma)}{\Lambda_- (\gamma)}
\ee
where $\Lambda_{\pm}$ are given by \eqref{consteig} and
$\gamma= \left(\E((1+x_1)^{-1})\right)^{-1}\in [3/2,2]$.
\een
\et
\bpr
We have, using \eqref{form}
\be\label{expr}
\frac{Y_n}{2}= \frac1n \log\left(\sum_{\alpha\in \{0,1\}^n}\left(\prod_{i=1}^n y_i^{\alpha_i-1} \right)
\tr (\caE(\alpha))\right)=:\psi(y_1,\ldots,y_n).
\ee
Using that the weights
\[
w(\alpha)=\frac{\prod_{i=1}^ny_i^{\alpha_i-1}\tr(\caE(\alpha))}
{\sum_{\alpha\in \{0,1\}^n}\biggl(\prod_{i=1}^n y_i^{\alpha_i-1}\biggr)\tr(\caE(\alpha))}
\]
are non-negative and sum up to one,
we compute the variation
\beq
L_i\psi &=& \sup_{y_1,\ldots,b_i,y_n}| \psi(y_1,\ldots,b_i,\ldots,y_n)-\psi(y_1,\ldots,y_n)|
\nonumber\\
&\leq &
\frac1n
\left(
\log\sum_{\alpha\in \{0,1\}^n} w(\alpha) \biggl(\frac{b_i}{y_i}\biggr)^{\alpha_i-1}
\right)\nonumber\\
&\leq & \frac1n \sup_{a,b\in [1/2,1]} \log \biggl|\frac{b}{a} \biggr| = \frac{\log 2}{n}.
\eeq
Statement 1 is then an application of the
Azuma-Hoeffding inequality.
The first part of Statement 2 follows from Oseledec's ergodic Theorem \cite{os}, together
with
\eqref{detM} and the law of large numbers.
To prove Statement 3, start from \eqref{form}, \eqref{3}, and use Jensen's inequality, combined with
the mutual independence of $y_i=1+x_i,1\le i\le n$, to obtain
\beq
&&\lim_{n\to\infty}\E \Bigl(\frac1{2n} \log \left(\frac{\lambda_+(\caM)}{\lambda_-(\caM)}\right) \Bigr)
\nonumber\\
&=&
\lim_{n\to\infty}\E\left(\frac1n \log\sum_{\alpha\in \{0,1\}^n}\left(\prod_{i=1}^n y_i^{\alpha_i-1} \right)
\tr(\caE(\alpha))\right)\nonumber\\
&\leq &
\lim_{n\to\infty}\frac1n\log\left(\sum_{\alpha\in \{0,1\}^n}\left(\prod_{i=1}^n \E\left(y_i^{\alpha_i-1}\right)\right)
\tr(\caE(\alpha))\right)\nonumber\\
&=&\lim_{n\to\infty}
\frac1n\log\left(\sum_{\alpha\in \{0,1\}^n}\left(\prod_{i=1}^n \left(\E(y_1^{-1})\right)^{(-1)(\alpha_i-1)}\right)
\tr(\caE(\alpha))\right)
\nonumber\\
&=&
\frac12\log\left(\frac{\Lambda_{+}(\gamma)}{\Lambda_-(\gamma)}\right)
\eeq
with $\gamma^{-1}= \E(y_1^{-1})$, and $\Lambda_{\pm}(\gamma)$ given by \eqref{consteig}.
Because of \eqref{detM},
 we finally derive \eqref{prodeig} by the law
of large numbers.
\epr
\section{Covariance and avalanche sizes}\label{covavsize}
\subsection{Quenched and annealed covariance}\label{quancov}
\bt\label{unifcov}
Let $\caT$ be a finite or infinite subtree of the full binary tree
$\caB_\infty$. As before, denote by $\mu_{\caT}$ the uniform
measure on the recurrent configurations $\caR_{\caT}$  of the sandpile
model on $\caT$ and let $u,v\in \caT$ be at mutual distance $n$.
Then we have the following estimate for the
two point correlation function
\be\label{covesti}
|Cov(u,v,\caT)|:=\left|\int\eta_{u} \eta_v \mu_\caT (d\eta)-
\int \eta_u \mu_\caT (d\eta)\int \eta_v\mu_{\caT} (d\eta)\right|
\leq
C \gamma^n.
\ee
In particular, $Cov(u,v,\caT)$ is absolutely summable, uniformly in the choice
of the subtree $\caT$:
\be\label{unifesti}
\sup_{\caT\subset\caB_\infty}\sum_{v\in \caT} |Cov(u,v,\caT)| <\infty.
\ee
\et
\bpr
The first statement \eqref{covesti} follows from
the expression of the covariance in Subsection \ref{transmatapp}, eq. \eqref{covbinary}, together with Proposition \ref{propo}.
Since the number of vertices at distance $n$ in $\caT$ from the origin is bounded
from above by $2^n$, and $\gamma<1/2$, \eqref{unifesti} follows
from \eqref{covesti}.
\epr
\bt\label{annealcov}
Let $\caT[o,n]$ be the stationary binary branching process with
reproduction probability $p$, conditioned to have a path
from its root $o$ to a vertex at distance $n$. Let $\mu_{\caT[o,n]}$ denote
the uniform measure on recurrent configurations on $\caT[o,n]$.
Let $Cov(o,v(n),\caT[o,n])$ denote the covariance of the height
variables at $o$ and at a vertex $v(n)$ at distance $n$ from the root $o$.
Then we have the following annealed lower bound on the
covariance.
\be\label{annealedcov}
\lim_{n\to\infty}\frac1n \left(-\log \E \left( |Cov(o,v(n),\caT[o,n])\right)| \right)\leq \log
\left(\frac{\Lambda_+ (\gamma)}{\Lambda_-(\gamma)}\right)
\leq \left(\frac{4+\sqrt{7}}{4-\sqrt{7}}\right),
\ee
where $\gamma^{-1}=\E(y_1^{-1})$ and $y_1=1+x_1$, with $x_1$  distributed
according to the measure of Proposition \ref{contractionbinary}, and $\Lambda_{\pm}$ given
in \eqref{consteig}.
\et
\bpr
The result follows from
the expression of the covariance in Subsection \ref{transmatapp}, eq. \eqref{covbinary},
together with \eqref{annealup}.
\epr
\subsection{Avalanche sizes}\label{avsiz}
We start by defining a matrix associated to an avalanche
cluster. Roughly speaking, the probability that
the avalanche coincides with the cluster is the
inverse of the maximal eigenvalue of this matrix.
 Let $\caT\subset\caB_\infty$ be a finite or infinite
subtree of the full rootless binary tree, containing the origin $o$.
Let $\caC$ be a connected subset of $\caT$ containing the origin.
\bd\label{matrdef}
Let $|\caC|=n$.
The matrix associated to $\caC$ in $\caT$, denoted by $M(\caC)$
is defined as follows. To the vertices in $\caC$ are
associated $n+2$ subtrees
$T_1,\ldots,T_{n+2}$ with corresponding characteristic ratios
$x(T_i)$.
Then
\be\label{matrixc}
M(\caC)= \prod_{i=1}^{n+2}
\left(
\begin{array}{cc}
1+x(T_i) & 1+x(T_i)\\
1 & 2+x(T_i)
\end{array}\right).
\ee
\ed
\bl\label{avlem}
As before, for a stable height configuration
$\eta$, let $Av(o,\eta)$ denote the avalanche
caused by addition of a single grain at the origin in $\eta$, and
$\mu_{\caT}$ the uniform measure on recurrent configurations $\caR_{\caT}$
on $\caT$.
Then there exist constants $c_1,c_2>0$ such that
\be
c_1 \left(\lambda_+ (M(\caC))\right)^{-1}\leq\mu_\caT (Av(o,\eta)=\caC)\leq c_2 (\lambda_+
\left( M(\caC)\right))^{-1}
\ee
where $M(\caC)$ is the matrix of \eqref{matrixc} associated to $\caC$ of Definition
\ref{matrdef}.
\el
\bpr
It follows from the expression for $\mu_\caT (Av(o,\eta)=\caC)$
in Subsection \ref{transmatapp}, eq. \eqref{avdis}.
\epr
\bd
We denote by $A_n(o,\caT)$
the number of connected subsets of edges of $\caT$
containing the origin $o$ and of cardinality $n$.
\ben
\item
The growth rate is  defined
as
\be\label{grow}
\kappa (o,\caT)=\limsup_{n\to\infty}\frac1n\log A_n(o,\caT).
\ee
\item We define the averaged growth rate as
\[
\bar{\kappa}=\limsup_{n\to\infty}\frac1n\log\E(A_n(o,\caT)).
\]
\een
\ed
The growth rate gives the dominant exponential
factor in the number $A_n(o,\caT)$. E.g. if $\caT$ is
the full binary tree, $\kappa=\log 4$, since
\be\label{fulltree}
A_n(o,\caB_\infty)= C4^n n^{-3/2}(1+o(1))
\ee
see also \cite{dhar}, Subsection \ref{cr+rec}.
For a stationary binary branching process we have the
upper bound
\be\label{branchup}
\E(A_n(o,\caT))\leq C p^n A_n(o,\caB_\infty).
\ee
For an exact expression $\E(A_n(o,\caT))$ from which inequality \eqref{branchup}
follows immediately, we refer to the Appendix, Proposition \ref{prop:7.2}.
\bt\label{decaythm}
\ben
\item  There exists $C>0$ such that for all $n\geq 1$ we have
\be\label{decay1}
\mu_\caT (|Av(o,\eta)|=n)\leq CA_n(o,\caT) 4^{-n} 2^{n(16/25)}.
\ee
In particular if the tree $\caT$ has growth rate $\kappa(\caT)< \frac{34}{25}\log 2$, 
then the avalanche size decays exponentially.
\item For the stationary binary branching tree with branching probability $p$,
we have the estimate
\be\label{decay2}
\E (\mu_\caT (|Av(o,\eta)|=n))\leq C 2^{n(16/25)} \biggl( \frac{p+\sqrt{p}}{2}\biggr)^{n}.
\ee
In particular, for
\[
p< 0.54511...
\]
the averaged (over the
realization of the tree) probability of an avalanche size $n$
decays exponentially in $n$.
\item For the binomial branching tree  we have
\be
\E (\mu_\caT (|Av(o,\eta)|=n))\leq C 2^{-n(16/25)}p^n,
\ee
\een
In particular, for
\[
p< 2^{-16/25} \approx 0.641713...
\]
the averaged (over the
realization of the tree) probability of an avalanche size $n$
decays exponentially in $n$.
\et
\bpr
The first result, i.e., \eqref{decay1},
follows from Lemma \ref{avlem}, and the estimate \eqref{lowertrace} which
gives the following lower bound on the largest eigenvalue
\be\label{lam+}
\lambda_+(M(\caC))= \tr (M(\caC)(1+o(1))\geq C4^n 2^{-n(16/25)}.
\ee
The second result follows from Lemma \ref{avlem},
\eqref{lam+}, \eqref{branchup},  \eqref{fulltree} and Proposition \ref{largen} below.
The solution of
\be
\frac{p+\sqrt{p}}{2 } = 2^{-16/25} \label{functP}
\ee
is given by
\be
p \approx 0.54511
\ee
The l.h.s. of \eqref{functP} as a function of $p$ is monotone, which yields that for $p < 0.54511$
the expected number of avalanches of size $n$ decays exponentially.

The third statement follows from
Lemma \ref{avlem},
\eqref{lam+}, \eqref{branchup},  \eqref{fulltree} and Remark \ref{remarkbinomial}.
\epr
\br
{}From Theorem \ref{decaythm} we conclude that 
avalanche sizes decay exponentially for $p$ small enough. We know that on the full
binary tree (corresponding to $p=1$) avalanche sizes have
a power law decay. It is a natural question whether the transition
between exponential and power law decay occurs at some
unique non-trivial value $p_c\in (1/2, 1)$.
 This question is related to the behavior of a a random walk on the full binary tree, killed upon exiting a random subtree.
Would this random walk have a survival time with a finite exponential moment, then we would have $p_c=1$. Notice that for the corresponding
problem on the lattice $\mathbb{Z}^d$, this random walk {\em does not have a finite exponential moment} of its survival time, because the tail of the
survival time is stretched exponential (Donsker-Varadhan tail $e^{-c t^{d/(d+2)}}$). However, on the tree in the annealed case corresponding to Theorem
\ref{decaythm}, the decay of the tail of the survival time is not straightforward. 
Formally taking the limit $d \rightarrow \infty$ in the Donsker-Varadhan tail suggests that 
on the tree, the survival time of this random walk has a finite exponential moment, which points into the
direction $p_c=1$. See also \cite{chen} for trapped random walk on a tree. 
\er
\noindent{\bf Acknowledgements}.
We would like to thank Remco van der Hofstad for useful discussions. F.R. and W.M.R. thank NWO for financial support in the project ``sandpile models in neuroscience'' within the STAR cluster.
This work was supported by
ANR-2010-BLAN-0108 and by  NWO. We also thank laboratoire MAP5
at Universit\'e Paris Descartes, Universities of Leiden and
Delft for financial support and hospitality.

\section{Appendix: The expected number of clusters containing the origin in a random binary tree}\label{appendix}
\bp\label{largen}
Let $p$ be the probability that a vertex has 2 children. Furthermore let $A_n(o, \caT^{o})$ denote the number of
connected clusters of $n$ edges containing $o$, given a realization of a tree $\caT^{o}$
with root $o$. Then there exists a constant $C>0$ such that for $n\geq 1$ we have
\begin{equation}\label{eq:largen}
\mathbb{E}(A_n(o,\caT^{o})) \leq C (n+1) 4^n \biggl (\frac{p+\sqrt{p}}{2} \biggr )^n.
\end{equation}
\ep
\br\label{remarkbinomial}
For the binomial branching tree, where every vertex has 2 children 
with probability $p^2$ and 1 child with probability $2p(1-p)$ 
(see Remark \ref{rk:intro-binomial}), we have the upper bound
\begin{equation}
\begin{split}
\mathbb{E}(A_n(o,\caT^{o})) \leq C4^n p^n
\end{split}
\end{equation}
which follows from the observation that the number of $n$ vertex animals containing the root
of the full binary tree is bounded from above by $4^n$, and in the case of binary branching tree, each
of these vertices is present with probability $p$ independently, which gives the factor $p^n$.
\er

\bpr
If we define $a_n$ to be the expected number of connected clusters containing $o$ and with $n$ vertices, then
$a_n$ satisfies the recursion relation
\be
 a_n =  \1_{\lbrace n= 1 \rbrace} + p \1_{\lbrace n \geq 2 \rbrace}\sum_{i=0}^{n-1}a_ia_{n-i-1},\qquad n\ge 1. \label{an}
\ee
Furthermore, by definition we put $a_0=1$. Thus, going from vertices to edges,
\be
\E(A_n(o,\caT^{o})) = a_{n+1}. \label{rel}
\ee
Introduce then the associated generating function:
\be
\begin{split}
A(x)=\sum_{n=0}^\infty a_n x^n & = 1 + x + p\sum_{n=2}^{\infty} \sum_{i=0}^{n-1}a_i a_{n-i-1} x^n \\
& = 1 + x + px(A^2(x)-1).\label{ax2}
\end{split}
\ee
For $\displaystyle{x\in\left[0,\frac{-p+\sqrt{p}}{2p(1-p)}\right]}$, this power series is convergent and, using
$A(0)=1$, it is given by
\be
A(x) = \frac{1}{2px} \biggl (1 - \sqrt{1-4px(1+x(1-p))} \biggr ),\label{ax3}
\ee
for $\displaystyle{x\in\left (0,\frac{-p+\sqrt{p}}{2p(1-p)}\right]}$.
Use (for $z$ such that $4z<1$)
\be
\begin{split}
\sqrt{1-4z} & = 1 - 2\sum_{n=1}^{\infty}{2n-2 \choose n-1} \frac{1}{n}z^n
\end{split}
\ee
to obtain
\be
\begin{split}
A(x) & = \sum_{n=0}^{\infty} \sum_{j=0}^{n+1} {2n \choose n} {n+1 \choose j} (1-p)^j \frac{p^n}{n+1}x^{n+j}.
\end{split}
\ee
We put $k=n+j$, then $n=k-j$ and $j \leq n+1$ yields $j \leq \lfloor (k+1)/2\rfloor$, 
hence the expected number of clusters $a_k$ containing the origin is equal to
\be
a_k = p^k\sum_{j=0}^{\lfloor (k+1)/2 \rfloor} b_{j,k} \label{summands}
\ee
where
\beq\nonumber
b_{j,k}&=&{2 (k-j) \choose k-j} {k-j+1 \choose j} \biggl(\frac{1-p}{p} \biggr)^j \frac{1}{k-j+1}\\
&  =& \frac{(2(k-j)!)}{(k-j)!j!(k-2j+1)!} c^j\label{st}
\eeq
with $c= (1-p)/p$.
All terms in the sum \eqref{summands}
are exponentially large in $k$, hence the exponential growth of the sum is determined by its maximal term.
Using Stirling's approximation $n! \approx  n^n e^{-n}\sqrt{2\pi n}$ for the right hand side of \eqref{st},
we obtain the upper bound
\be\label{bkj}
b_{j,k}\leq  \frac{e}{\pi\sqrt{2}} c^j 4^{k-j} \frac{(k-j)^{k-j}}{j^{j}(k-2j)^{k-2j}}
\ee
Define, for $x\in [0,(k+1)/2]$,
\[
\phi(x)= \frac{e}{\pi\sqrt{2}} c^{x} 4^{k-x} \frac{(k-x)^{k-x}}{x^{x}(k-2x)^{k-2x}}.
\]
This function $\phi$ attains its maximum at $x= k(1-\sqrt{p})/2$, which combined with \eqref{bkj}, implies
\be\label{fj}
b_{j,k}\leq \frac{e}{\pi\sqrt{2}} 4^k \biggl (\frac{1+\sqrt{p}}{2\sqrt{p}} \biggr )^k.
\ee
Plugging \eqref{fj} into \eqref{summands} and bounding from above induces \eqref{eq:largen}.
\epr

In the following proposition we give an exact formula for $\mathbb{E}(A_n(o,\caT^{o}))$.
\bp\label{prop:7.2}
We have the following identity,
\be
\begin{split}
& \mathbb{E}(A_n(o,\caT^{o}))\\
& = \frac{p^{n+2}}{n+2} {2(n+1) \choose n+1}   {_2F_1}\biggl(-\frac{n+2}{2},-\frac{n+1}{2},\frac{1}{2}-(n+1),\frac{-(1-p)}{p} \biggr) \label{identity}
\end{split}
\ee
where $_2F_1(\cdot,\cdot,\cdot,\cdot)$ denotes the hypergeometric function defined as
\be
_2F_1(a,b,c,z) = \sum_{n=0}^{\infty} \frac{(a)_n(b)_n}{(c)_n}\frac{z^n}{n!}\label{hyp}
\ee
and where the Pochhammer symbol $(a)_n$ is defined by
\be
(a)_n = \begin{cases} 1 & \text{ if }  n=0, \\ a(a+1)...(a+n-1) & \text{ if } n> 0. \end{cases} \label{Poch}
\ee
\ep
\bpr
Let us first remark that for $|z| < 1$, real $a,b$ and $c \neq -m$ where 
$m \in \mathbb{N}$, the series in \eqref{hyp} is
well defined. In our case $z=-(1-p)/p$, $a=-(n+1)/2$, $b=1/2-(n+1)$ and 
$c = (n+1)-1/2$ ( see also  \cite{andrews} for more details about the
hypergeometric function).
The claim \eqref{identity} follows from \eqref{rel} and \eqref{summands} 
once we show the identity
\be
\begin{split}
&\sum_{j=0}^{\lfloor (k+1)/2 \rfloor} {2 (k-j) \choose n+1-j} {k-j+1 \choose j} \biggl(\frac{1-p}{p} \biggr)^j
\frac{1}{k-j+1}  \\
& = \frac{1}{k+1} {2k \choose k}   {_2F_1}\biggl(-\frac{k+1}{2},-\frac{k}{2},\frac{1}{2}-k,\frac{-(1-p)}{p} \biggr). \label{identity1}
\end{split}
\ee
If $a$ or $b$ are negative integers, using \eqref{Poch} we see 
that the series \eqref{hyp} defining ${_2F_1}$ is a finite sum. Therefore,
\be
\begin{split}
 {_2F_1}\biggl(-\frac{k+1}{2},-\frac{k}{2},\frac{1}{2}-k,-\frac{1-p}{p} \biggr) &= \sum_{j=0}^{\lfloor k+1/2 \rfloor} \frac{\left(-\frac{k+1}{2} \right)_j \left(-\frac{k}{2} \right)_j}{\left(\frac12 - k\right)_j} \frac{(-1)^j}{j!}\frac{(1-p)^j}{p^j}.
\end{split}
\ee
We have
\be
(a)_n=\frac{\Gamma(a+n)}{\Gamma(a)}.
\ee
Next we use the functional identities for the Gamma function (see \cite{temme}),
\be
\Gamma(z)\Gamma(-z) = \frac{-\pi}{z\sin(\pi z)},
\ee
and recurrence relation
\be
z\Gamma(z) = \Gamma(z+1)\label{rec}
\ee
to rewrite
\be
\frac{\left(-\frac{k+1}{2}\right)_j \left(-\frac{k}{2}\right)_j}{ \left(\frac12 - k \right)_j} \label{anP}
\ee
in terms of Gamma functions with positive arguments. Thus
\be
\begin{split}
\left(-\frac{k+1}{2}\right)_j & = \frac{\Gamma\left(-\frac{k+1}{2} + j\right)}{\Gamma\left(-\frac{k+1}{2}\right)} \\
& = \frac{\Gamma\left(\frac{k+1}{2} + 1\right)\sin\left(\pi \frac{k+1}{2}\right)}{\Gamma\left(\frac{k+1}{2} + 1 - j\right)
\sin\left(\pi \left(\frac{k+1}{2} - j\right)\right)} \\
\end{split}
\ee
which gives
\be
\frac{\left(-\frac{k+1}{2}\right)_j \left(-\frac{k}{2}\right)_j}{ \left(\frac12 - k \right)_j} = A \times B
\ee
with
\be
A = \frac{\Gamma\left(\frac{k+1}{2} + 1\right)\Gamma\left(\frac{k}{2}+1\right) \Gamma\left(k -j +\frac12\right)}
{\Gamma\left(\frac{k+1}{2} + 1 - j\right)\Gamma\left(\frac{k}{2}+1-j\right)\Gamma\left(k+\frac{1}{2}\right)}
\ee
and
\be
B = \frac{\sin\left(\pi\left(\frac{k+1}{2}\right)\right)\sin\left(\frac{\pi k}{2}\right)
\sin\left(\pi\left(k-\frac12 -j\right)\right)}{\sin\left(\pi\left(\frac{k+1}{2}  - j\right)\right)
\sin\left(\pi\left(\frac{k}{2}-j\right)\right)\sin\left(\pi\left(k-\frac12\right)\right)}.
\ee
Furthermore expanding numerator and denominator of this expression gives
$B=(-1)^j$.
To rewrite $A$, we use \eqref{rec} and the duplication formula
\be
\Gamma(z)\Gamma\left(z+ \frac12 \right) = 2^{1-2z}\sqrt{\pi}\Gamma(2z)
\ee
once for $z=(k+1)/2$:
\be
\begin{split}
\Gamma\left( \frac{k+1}{2} + 1 \right) \Gamma\left(\frac{k+1}{2} + \frac12 \right)  & = \sqrt{\pi}\frac{k+1}{2} 2^{-k} \Gamma(k+1)\\
& = \sqrt{\pi}2^{-k-1}\Gamma(k+2),
\end{split}
\ee
and another time for $z=(k+1)/2-j$:
\be
\begin{split}
\Gamma\left( \frac{k+1}{2} - j + 1\right)\Gamma \left( \frac{k}{2} - j + 1\right)  &= \sqrt{\pi} \biggl(\frac{k+1}{2} - j \biggr)2^{-k + 2j} \Gamma(k -2j +1), \\
& = \sqrt{\pi} 2^{-k + 2j -1} \Gamma(k-2j + 2).
\end{split}
\ee
Hence
\be
A = \frac{\Gamma(k+2)}{\Gamma\left(k+\frac12 \right)} \times 4^{-j} \times \frac{\Gamma \left(k-j + \frac12 \right)}{\Gamma(k-j+2)} \times \frac{\Gamma(k-j+2)}{\Gamma(k-2j+2)}
\ee
and since the Catalan numbers satisfy
\be
\frac{1}{k+1} {2k \choose k} = \frac{4^k}{\sqrt{\pi}} \frac{\Gamma\left(k+\frac12 \right)}{\Gamma(k+2)}
\ee
we can write
\be
\begin{split}
\frac{4^k}{\sqrt{\pi}} \frac{\Gamma \left(k+ \frac12 \right)}{\Gamma(k+2)}  & \times A \times B  \times \frac{(-1)^j}{j!} = \\
& \frac{4^{k-j}}{\sqrt{\pi}} \frac{\Gamma \left(k-j + \frac12 \right)}{\Gamma(k-j+2)} \times \frac{\Gamma(k-j+2)}{\Gamma(k-2j+2) \Gamma(j+1)}
\end{split}
\ee
which is equal to
\be
{2(k-j) \choose k-j} \times {k-j+1 \choose j} \frac{1}{k-j+1}
\ee
and yields the claim.
\epr
\br
For $p=1$, we recover the classical formula for the number of animals of $n$ edges containing the origin
in a binary tree (cf. Subsection \ref{cr+rec}):
\be
A_n(o,\caT^{o}) = \frac{1}{n+2}{2(n+1) \choose n+1} \underset{n \text{  large }} \approx C 4^n n^{-3/2}
\ee
for some constant $C$.
\er
\bpr
This follows from the fact that for $p=1$, $_2F_1(a,b,c,0)=1$.
\epr

\begin{thebibliography}{99}
%
\bibitem{andrews}
G.E. Andrews, R. Askey and R. Ranjan, Special functions, \textit{Encyclopedia of Mathematics and its Applications}, {\bf 71}, Cambridge University Press. ISBN 978-0-521-62321-6 (1999).
%
\bibitem{chen}
M. Chen, S. Yan and X. Zhou, The range of random walk on trees and related trapping problem, \textit{Acta Math. Applic. Sinica}
{\bf 13} (1997), 1--16.
%
\bibitem{dharo}
D. Dhar, Theoretical studies of self-organized criticality, \textit{Phys. A} {\bf 369} (2006),  29--70.
%
\bibitem{dharb}
D. Dhar, Self-organized critical state of sandpile automaton models, \textit{Phys. Rev. Lett.} {\bf 64} (1990), 1613--1616.
%
\bibitem{dhar}
D. Dhar and S.N. Majumdar, Abelian sandpile model on the Bethe lattice, \textit{J. Phys. A: Math. Gen.} {\bf 23}  (1990), 4333--4350.
%
\bibitem{dudley}
R.M. Dudley, Real analysis and probability, \textit{Cambridge University Press}, (2002).
%
\bibitem{cactus2}
G. Gauthier, Avalanche dynamics of the Abelian sandpile model on the expanded cactus graph,
preprint (2011), available at arxiv.org. 1110.6263v2.
%
\bibitem{smallworld}
K.I. Goh, D.S. Lee, B. Kahng, and D. Kim, Sandpile on Scale-Free Networks,
\textit{Phys. Rev. Lett.} {\bf 91} (2003), 148701.
%
\bibitem{jaraio}
A.A. J\'{a}rai, Thermodynamic limit of the abelian sandpile model on $\Z^d$. \textit{Mark. Proc. Rel. Fields} {\bf 11} (2005), 313--336.
%
\bibitem{lev}
A. Levina, J.M. Herrmann, and T. Geisel, Dynamical synapses causing
self-organized criticality in neural networks, \textit{Nature Physics}, {\bf 3} (2007), 857--860.
%
\bibitem{MRS}
C. Maes, F. Redig and E. Saada, The Abelian sandpile model on an infinite tree, \textit{Ann. Probab.} {\bf 30} (2002), 2081--2107.
%
\bibitem{MRSindien}
C. Maes, F. Redig and E. Saada, Abelian sandpile models in infinite volume.
{\sl Sankhya, the Indian Journal of Statistics}, {\bf 67}, no. 4 (2005), 634--661.
%
\bibitem{cactus}
M. Matter and T. Nagnibeda, Abelian Sandpile Model on Randomly Rooted Graphs and
Self-Similar Groups, preprint (2010), available at arxiv.org. 1105.4036.
%
\bibitem{os}
V. L. Osedelec, A Multiplicative Ergodic Theorem; Lyapunov Characteristic Number
for Dynamical Systems, \textit{Trans. Moscow Math. Soc.} {\bf 19} (1968), 197--231.
%
\bibitem{redigo}
F. Redig, Mathematical aspects of the abelian sandpile model,
\textit{Mathematical statistical physics, Les Houches Summer School}, Elsevier B. V., Amsterdam (2006), pp.\ 657--729.
%
\bibitem{temme}
N.M. Temme, Special Functions: An Introduction to the Classical Functions of Mathematical Physics, \textit{John Wiley \& Sons}, New York,  ISBN 0-471-11313-1 (1996).
%
\end{thebibliography}
\end{document}